\documentclass[review]{elsarticle}
\usepackage{amsmath}
\usepackage{picinpar,graphicx}
\usepackage{lineno,hyperref}

\modulolinenumbers[5]

\journal{Journal of \LaTeX\ Templates}

\newcommand{\bqa}{\begin{eqnarray}}
\newcommand{\eqa}{\end{eqnarray}}
\newcommand{\bqn}{\begin{eqnarray*}}
\newcommand{\eqn}{\end{eqnarray*}}
\newcommand{\be}{\begin{equation}}
\newcommand{\ee}{\end{equation}}
\newcommand{\md}{\mbox{d}}
\newcommand{\non}{\nonumber\\}

\newtheorem{thm}{Theorem}[section]
\newtheorem{lem}{Lemma}
\newtheorem{rmk}{Remark}[section]

\newtheorem{ass}{Assumption}

\newcommand{\um}{\underline{m}}

\begin{document}

\begin{frontmatter}

\title{The limits of the sample spiked eigenvalues for a high-dimensional generalized Fisher matrix and its applications}

\author[mymainaddress]{Dandan Jiang\fnref{myfootnote1}}
\fntext[myfootnote1]{Supported by NSFC (No.  11971371).}
\ead{jddpjy@163.com}

\author[mysecondaryaddress]{Jiang Hu\fnref{myfootnote2}}
\fntext[myfootnote2]{Supported by
NSFC (No. 11771073) and Foundation of Jilin Educational Committee (No. JJKH20190288KJ).}

\ead{huj156@nenu.edu.cn}

\author[mysecondaryaddress]{Zhiqiang Hou\corref{mycorrespondingauthor}}
\cortext[mycorrespondingauthor]{Cortrresponding author}

\ead{houzq399@nenu.edu.cn}

\address[mymainaddress]{School of Mathematics and Statistics,
	Xi'an Jiaotong University,
	No.28, Xianning West Road,
	Xi'an {\rm 710049}, China.}
\address[mysecondaryaddress]{School of Mathematics and Statistics, Northeast Normal University, 5268 Renmin Street, Changchun {\rm 130024}, China}

\begin{abstract}
A generalized spiked Fisher matrix is considered in this paper. We establish a criterion for the description of the support of the limiting spectral distribution of high-dimensional generalized Fisher matrix and study the almost sure limits of the sample spiked eigenvalues where the population covariance matrices are arbitrary which successively removed an unrealistic condition posed in the previous works, that is, the covariance matrices are assumed to be diagonal or diagonal block-wise structure. In addition, we also give a consistent estimator of the population spiked eigenvalues. A series of simulations are conducted that support the theoretical results and illustrate the accuracy of our estimators.
\end{abstract}


\end{frontmatter}


\section{Introduction} \label{Int}

 Consider the  spiked model involved with  two sample covariance matrices,
\be
{\boldsymbol\Sigma}_1={\boldsymbol\Sigma}_2+{\boldsymbol\Delta},\label{model}
\ee
where ${\boldsymbol\Sigma}_1$ and ${\boldsymbol\Sigma}_2$ are general covariance matrices and ${\boldsymbol\Delta}$ is a finite rank matrix.
This two-sample spiked model  has wide applications to many fields, including signal processing, regression analysis, etc.
To illustrate, we enumerate several basic problems, such as
testing the presence of signals and testing the number of signals in signal processing. Additionally, the Lawley-Hotelling trace criterion, the Bartlett-Nanda-Pillai trace criterion and the Roy Maximum root criterion are used in testing the linear regression hypothesis.
Under the alternative hypothesis,
these tests are based on the sample spiked eigenvalues of the Fisher matrix ${\boldsymbol S}_1{\boldsymbol S}_2^{-1}$, where ${\boldsymbol S}_1, {\boldsymbol S}_2$ are the sample covariance matrices corresponding to ${\boldsymbol\Sigma}_1$ and ${\boldsymbol\Sigma}_2$, respectively.
However, the sample spiked eigenvalues do not converge to their corresponding population spiked eigenvalues if the dimensionality $p$ goes to infinity.
Therefore, traditional testing methods and their asymptotic laws lose efficiency in such a case.
Thus, a study of the limits of sample spiked eigenvalues is necessary.

There are many works that investigate the spiked model in a high-dimensional setting.
As is well known, the  spiked model,  first proposed by \cite{Johnstone2001}, 
can be seen as a special case of ${\boldsymbol\Sigma}_2={\boldsymbol I}$ in (\ref{model}), which  has  the same approach as that of principal component analysis (PCA).
Then, some relevant works are devoted to improving the study of the one-sample spiked model, such as
\cite{BaiNg2002}, \cite{Baik2005}, \cite{BaikSilverstein2006}, \cite{Paul2007}, \cite{BaiYao2008,BaiYao2012},  \cite{Onatski2009,Onatski2012}, \cite{FanWang2015}, \cite{Caietal2017}.
Some related studies are also  devoted to investigations of PCA or FA, which can be seen as another way of understanding the spiked model.
Examples include \cite{HoyleRattray2004}, \cite{Nadler2008},
\cite{JungMarron2009}, \cite{Shen2013},  \cite{BerthetRigollet2013}, \cite{Birnbaum2013}, etc.
Recently,  \cite{JiangBai2019}
extended the work to a general case
and gave the limits and CLT for the sample spiked eigenvalues of a generalized covariance matrix.

In contrast, there are only a few studies related to the two-sample spiked model.
\cite{WangYao2017} assumed that ${\boldsymbol\Sigma}_1{\boldsymbol\Sigma}_2^{-1}$ is an identity matrix with a rank M perturbation or diagonal block independent and presented the limits of the  extreme eigenvalues of a high-dimensional spiked Fisher matrix.
In addition, \cite{JohnstoneOnatski2015}
described the  relationship between the two-sample spiked models with some classical statistical problems that lead to each of James' five cases in \cite{James1964}. In the alternative hypothesis, they focused on the two-sample spiked model with ${\boldsymbol\Delta}$ being a rank-one matrix that is used to derive the asymptotic power for testing the presence of a spike.
However, these works are all based on the simplified structure of the Fisher matrix and are limited in practice. First, the diagonal or diagonal blockwise assumption is an impractical assumption, which means that the spiked and non-spiked eigenvalues are generated from independent variables. Moreover, the rank-one assumption is the same as the fact that there is only one input signal. Thus, it is not applicable to other statistical inferences in signal processing, such as testing the number of signals.  Therefore, there is still room for improvement in these studies.

Note that the existing limiting laws for the spiked eigenvalues of the simplified Fisher matrix are established based on the normalized difference between the sample spiked eigenvalues and their limits. Thus, we extend to a generalized spiked Fisher matrix and focus on
the first step for the tests on the spikes, which is to calculate the limits of the sample spiked eigenvalues  with high dimensionality $p$.
As a natural consequence, the estimators of the population spikes are also obtained, which can be used to restore the concerned matrix structure. Moreover, the estimated population spikes can represent the strength of the input signals.

The main contributions of the paper include: established a criterion for the description of the support of the limiting spectral distribution of high-dimensional generalized Fisher matrix; established the almost sure limits of the sample spiked eigenvalues where the population covariance matrices are arbitrary which successively removed an unrealistic condition posed in the previous works, that is, the covariance matrices are assumed to be diagonal or diagonal block-wise structure. In addition, we also give a consistent estimator of the population spiked eigenvalues. A series of simulations are conducted that support the theoretical results and illustrate the accuracy of our estimators.

The paper is organized as follows. In Section~\ref{Pre}, we present the almost sure limits of the sample spiked eigenvalues for a high-dimensional generalized Fisher matrix and establish a criterion for the description of the support of the limiting spectral distribution of high-dimensional generalized Fisher matrix, which are the main results of the paper. Section~\ref{New} gives estimators of the population distant spiked eigenvalues for the generalized Fisher matrix. In Section~\ref{Sim}, we conduct simulations that support the theoretical results and illustrate the accuracy of the estimators of the population distant spiked eigenvalues. Technical lemmas and proofs are postponed to the Supplementary Material.

\section{The limits of the sample spiked eigenvalues for a Generalized Spiked Fisher matrix.
  \label{Pre}}
Assume that
\bqn
{\mathbf X}&=&({\mathbf x}_1,\cdots,{\mathbf x}_{n_1})=\left(x_{ij}\right), 1\leq i \leq p, ~ 1\leq j \leq n_1,\\
{\mathbf Y}&=&({\mathbf y}_1,\cdots,{\mathbf y}_{n_2})=\left(y_{kl}\right), 1\leq k \leq p, ~ 1\leq l \leq n_2
\eqn
are two independent  $p$-dimensional arrays with  components having zero mean  and  identity variance.
Denote
${\boldsymbol\Sigma}_1^{1/2}{\mathbf X}$ and  ${\boldsymbol\Sigma}_2^{1/2}{\mathbf Y}$ as two independent  samples with  two population covariance matrices, where ${\boldsymbol\Sigma}_1$ and ${\boldsymbol\Sigma}_2$ are two general nonnegative definite matrices.
Let ${\mathbf T}_p={\boldsymbol\Sigma}_1^{1/2}{\boldsymbol\Sigma}_2^{-{1/2}}$ and further assume that the spiked eigenvalues of ${\mathbf T}_p^*{\mathbf T}_p$ are   scattered into spaces of a few  bulks with the largest allowed to tend to infinity.
Thus, for the corresponding  sample covariance matrices of the two observations, 
\begin{equation}{\mathbf S}_1
=\frac{1}{n_1}{\boldsymbol\Sigma}_1^{1\over 2}
{\mathbf X}{\mathbf X}^*{\boldsymbol\Sigma}_1^{1 \over 2}  \quad \text{and} \quad   {\mathbf S}_2
=\frac{1}{n_2}{\boldsymbol\Sigma}_2^{1 \over 2}
{\mathbf Y}{\mathbf Y}^*{\boldsymbol\Sigma}_2^{1 \over 2},\label{S12}
\end{equation}
the matrix ${\mathbf F}={\mathbf S}_1{\mathbf S}_2^{-1}$
is the so-called generalized Fisher matrix, where the condition $n_2>p$ is necessary for the invertible
matrix ${\mathbf S}_2$. Because the matrix ${\mathbf F}$ has the same nonzero eigenvalues as those of the matrix,
\begin{equation}
{\mathbf F}={\mathbf T}_p^{*}\widetilde {\mathbf S}_1{\mathbf T}_p\widetilde {\mathbf S}_2^{-1},\label{F}
\end{equation}
where $\widetilde {\mathbf S}_1=n_1^{-1}{\mathbf X}{\mathbf X}^*$ and $\widetilde {\mathbf S}_2={n_2}^{-1}{\mathbf Y}{\mathbf Y}^*$  are the standardized sample covariance matrices,
we investigate the Fisher matrix ${\mathbf F}={\mathbf T}_p^{*}\widetilde {\mathbf S}_1{\mathbf T}_p\widetilde {\mathbf S}_2^{-1}$ instead. If there is no confusion, we will still use the notation ${\mathbf F}$.

Furthermore,  we  assume that the spectrum of  ${\mathbf T}_p^*{\mathbf T}_p$ is listed in descending order
as below:
\begin{equation}
\beta_{p,1}, \cdots,  \beta_{p,j},\cdots,\beta_{p,p}.\label{array}
\end{equation}
Denote the spikes as $\beta_{p,j_k+1}= \cdots= \beta_{p, j_k+m_k} \overset{def}{=}\alpha_k$ with $j_k's$ being  arbitrary  ranks in the array
(\ref{array}); then,  the population spiked eigenvalues
$\alpha_1, \cdots, \alpha_K$ with multiplicity $m_k, k=1,\cdots,K$  are  aligned arbitrarily in groups among all the eigenvalues, satisfying $m_1+\cdots+m_K=M$, a fixed integer.
In addition,  the spiked eigenvalues are allowed  to be infinity.  Under these general assumptions, the matrix ${\mathbf F}$ is called a generalized spiked Fisher matrix.

To study the limiting behaviors 
of  the distant sample spiked eigenvalues of the generalized Fisher matrix ${\mathbf F}$,  some necessary assumptions are detailed as follows:

\begin{ass}
	\label{assA}
	Let $\{x_{ij}, i=1,\dots,p,j=1,\dots,n_1\}$ be a set of independent and identically distributed (i.i.d.) random variables with mean 0, variance 1 and finite fourth moments. Analogously, let $\{y_{ij}, i=1,\dots,p,j=1,\dots,n_2\}$ be another set of  i.i.d. random variables that are independent of $\{x_{ij}\}$ with mean 0, variance 1 and finite fourth moments. If they are   complex,
	$E x_{ij}^2=0$ and $E y_{ij}^2=0$ are required.
	\end{ass}
\begin{ass}
	\label{assB}
	The matrix ${\mathbf T}_p={\boldsymbol\Sigma}_1^{1 \over 2}{\boldsymbol\Sigma}_2^{-{1 \over 2}}$ is nonrandom and has all its eigenvalues bounded except for a fixed number of eigenvalues that are allowed to be infinite at a rate of $o(n^{1/6})$.
	Moreover,  the empirical spectral distribution of $\{{\mathbf T}_p^*{\mathbf T}_p\}$,  
	denoted by $H_n$, tends to proper probability measure $H$ as $\min( p,n_1,n_2)\rightarrow \infty$.
\end{ass}
\begin{ass}
	\label{assC}
	Assume that
	$c_{n_1}=p/n_1\rightarrow c_1\in (0, \infty)\quad \text{and} \quad  c_{n_2} =p/n_2 \rightarrow c_1\in (0,1)$ 
	as  $\min ( p,n_1,n_2)\rightarrow \infty$.
\end{ass}

Our first aim is to investigate the limits of the sample spiked eigenvalues associated with $\alpha_i$ for a high-dimensional generalized Fisher matrix. 
To be specific,  for any measure $\mu$ on $\mathcal{R}$, we denote  the support  of $\mu$ as $\mathcal{G}_{\mu}$, a closed set.
Then,
the eigenvalue $\beta_{p,j}$ is a spiked eigenvalue if $\beta_{p,j} \notin \mathcal{G}_{H}$, where $H$ is the limiting spectral distribution of  ${\mathbf T}_p^*{\mathbf T}_p$.
To avoid possible confusion when the eigenvalues vary with the dimensionality $p$, we define the eigenvalues $\beta_{p,j}$ satisfying
$d(\beta_{p,j}, \mathcal{G}_{H})> \delta$ as the spiked eigenvalues, where $d$ is a predefined distance function and $\delta$ is a preselected positive constant.

Let $\mathcal{J}_k$ be the set of  ranks  of  $\alpha_k$ with multiplicity $m_k$ among all the eigenvalues of ${\mathbf T}_p^*{\mathbf T}_p$, {\rm i.e.,}
\[\mathcal{J}_k=\{ j_k+1,\ldots, j_k+m_k\}.\]
The sample eigenvalues of the generalized spiked Fisher matrix ${\mathbf F}$ are  arranged in descending order as
\[\lambda_{p,1}({\mathbf F}), \cdots,  \lambda_{p,j}({\mathbf F}), \cdots,  \lambda_{p,p}({\mathbf F}).\]
Let
	\[\varrho_k =
	\left\{
	\begin{array}{cc}
	\psi(\alpha_k),& \mbox{ if } \psi'(\alpha_k)>0,   \\
	\psi(\underline\alpha_k), &   \mbox{ if there exists } \underline\alpha_k \mbox{ such that } \psi'(\underline\alpha_k)=0
	\\ &\mbox{ and }\psi'(t)<0,  \mbox{ for all } \alpha_k\le t<\underline{\alpha}_k
	\\
	\psi(\overline\alpha_k), &   \mbox{ if there exists } \overline\alpha_k \mbox{ such that }\psi'(\overline\alpha_k)=0\\
	& \mbox{ and }
	\psi'(s)<0, \mbox{ for all } \overline{\alpha}_k<s\le\alpha_k
	\end{array}
	\right.
	\]
	where
	\be
	\psi(\alpha_k)=\frac{\alpha_k\left(1-c_1\int \displaystyle\frac{t}{t-\alpha_k}{\rm d}H(t)\right)}{1+c_2\int\displaystyle\frac{\alpha_k}{t-\alpha_k}{\rm d}H(t)}.\label{psik}
	\ee
Then, for each spiked eigenvalue $\alpha_k$ with multiplicity $m_k, k=1,\cdots, K$ associated  with sample eigenvalues
$\{\lambda_{p,j}({\mathbf F}), j \in \mathcal{J}_k\}$, we have
the following theorem. The proof is postponed to the Supplementary Material.


\begin{thm}\label{F-limits}
	Under Assumptions 1-3, for any integer  $k \in \{1,\dots,K\}$
	and all $ j \in {\mathcal J}_k$, as  $\min ( p,n_1,n_2)\rightarrow \infty$, we have that  $\lambda_{p,j}/\varrho_k-1\to0$  almost surely.
\end{thm}

\begin{rmk} \label{rmk1}
	Theorem \ref{F-limits} presents the limits of the sample eigenvalues associated with the population spike eigenvalues ${\alpha}_k$, where the involved $H$ is a general distribution different from the existing results such as those in \cite{WangYao2017}.  Theorem 3.1 in \cite{WangYao2017} is a special case of Theorem \ref{F-limits} when the limiting spectral distribution $H$ degenerates to $\delta_{\{1\}}$ with
	\[\psi(\alpha_k)=\frac{\alpha_k\left(1-\alpha_k-c_1\right)}{1-\alpha_k+c_2\alpha_k}.\]
\end{rmk}

In Theorem \ref{F-limits}, the $\alpha_k$'s satisfying $\psi'(\alpha_k)>0$ are called distant spiked eigenvalues, and the other two cases are called close spiked eigenvalues. The following two theorems give a criterion for the description of the support of the limiting spectral distribution of high-dimensional generalized Fisher matrix, in other words, they provide the close relationship between the population spike eigenvalues $\alpha_k$ and the limits of the sample outlier eigenvalues associated with $\alpha_k$, and they can help us complete the proof of Theorem \ref{F-limits}. In fact, these results are independent from the previous results and should have their own interest. The details of the proof are deferred to Supplementary Material.

Let $\mathbf n=(n_1, n_2), \mathbf c=(c_1, c_2)$ and $F_{\mathbf n}$ be  the empirical spectral distribution of ${\mathbf F}$, which  converges to a limiting spectral distribution $F^{\mathbf c, H}$.  Denote   $\mathcal{G}_{F^{\mathbf c, H}}$ as  the supporting set of the LSD  $F^{\mathbf c, H}$
and  $\mathcal{G}^c_{F^{\mathbf c, H}}$ as its complement. Then, we have


\begin{thm}\label{thm2}
	If $\lambda \in \mathcal{G}^c_{F^{\mathbf c,H}}$, then there exists $\alpha$  such that $\lambda=\psi(\alpha)$ and
	\begin{center}
		\begin{enumerate}
			\item[{\rm (i)}] $\alpha \in \mathcal{G}^c_{H}$, and $\alpha \neq 0$ such that the $\psi$ in  (\ref{psik}) is well defined.
			\item[{\rm (ii)}] $1-c_2\int\displaystyle\frac{\alpha^2 {\rm d} H(t)}{(t-\alpha)^2}>0$,
			\item[{\rm (iii)}] $\psi'(\alpha)>0$.
		\end{enumerate}
	\end{center}
\end{thm}

\begin{thm}\label{thm3}
	If the following conditions hold, {\rm i.e.,}
	\begin{center}
		\begin{enumerate}
			\item[{\rm (i)}] $\alpha \in \mathcal{G}^c_{H}$, and $\alpha \neq 0$ such that the $\psi$ in  (\ref{psik}) is well defined.
			\item[{\rm (ii)}] $1-c_2\int\displaystyle\frac{\alpha^2 {\rm d} H(t)}{(t-\alpha)^2}>0$,
			\item[{\rm (iii)}] $\psi'(\alpha)>0$.
		\end{enumerate}
	\end{center}
	then $\lambda \in \mathcal{G}^c_{F^{\mathbf c,H}}$, where $\lambda=\psi(\alpha)$.
\end{thm}

%

\section{Estimators of the population distant spiked eigenvalues.} \label{New}
For the generalized Fisher matrix ${\mathbf F}={\mathbf T}_p^{*}\widetilde {{\mathbf S}}_1{\mathbf T}_p\widetilde {{\mathbf S}}_2^{-1}$ defined in (\ref{F}), denote
the singular value decomposition of ${\mathbf T}_p$ as
\be {\mathbf T}_p= {\mathbf U}\left(
\begin{array}{cc}
	{\mathbf D}_1^{1/ 2} & \mathbf 0    \\
	\mathbf 0 & {\mathbf D}_2^{1/2}
\end{array}
\right){\mathbf V}^{*},
\label{SVD}
\ee
where ${\mathbf U}, {\mathbf V}$ are  unitary (orthogonal for the real case) matrices,   ${\mathbf D}_1$ is a diagonal matrix of the $M$ spiked eigenvalues of  the generalized spiked Fisher matrix ${\mathbf F}$
and ${\mathbf D}_2$ is the diagonal matrix of the non-spiked eigenvalues  with bounded components. Consider the $k$th bulk of the sample spiked eigenvalues  of ${\mathbf F}$, $\lambda_{p,j}, j \in \mathcal{J}_k$,  which satisfy the following eigen-equation
\[0=|\lambda_{p,j}{\mathbf I}- {\mathbf F}|=\left|\lambda_{p,j}{\mathbf I}-{\mathbf V}{\rm diag}({\mathbf D}_1^{1/2},{\mathbf D}_2^{1/2})
{\mathbf U}^* \tilde {{\mathbf S}}_1{\mathbf U}
{\rm diag}({\mathbf D}_1^{1/2},{\mathbf D}_2^{1/2}){\mathbf V}^*\tilde {{\mathbf S}}_2^{-1}\right|.\]
Partition the  two matrices, ${\mathbf U}, {\mathbf V}$, in the way of the matrix ${\mathbf D}={\rm diag}({\mathbf D}_1^{1/2},{\mathbf D}_2^{1/2})$; then, it
is equivalent to
\begin{eqnarray*}
	0&=&|\lambda_{p,j}{\mathbf V}^*\tilde {{\mathbf S}}_2{\mathbf V}\!-\!{\rm diag}({\mathbf D}_1^{1/2},{\mathbf D}_2^{1/2}){\mathbf U}^* \tilde {{\mathbf S}}_1 {\mathbf U}{\rm diag}({\mathbf D}_1^{1/2},{\mathbf D}_2^{1/2})|\\
	&=&| \lambda_{p,j} {\mathbf V}_2^* \tilde {{\mathbf S}}_2 {\mathbf V}_2\!-\!{\mathbf D}_2^{1/2}{\mathbf U}_2^* \tilde {{\mathbf S}}_1{\mathbf U}_2{\mathbf D}_2^{1/2}||{\mathbf K}(\lambda_{p,j})|
\end{eqnarray*}
where
\begin{eqnarray}\label{Kn}
{\mathbf K}(\lambda_{p,j})&=&\lambda_{p,j} {\mathbf V}_1^*\tilde {{\mathbf S}}_2 {\mathbf V}_1\!-\!{\mathbf D}_1^{1/2}{\mathbf U}_1^* \tilde {{\mathbf S}}_1{\mathbf U}_1{\mathbf D}_1^{1/2}\!-\!(\lambda_{p,j} {\mathbf V}_1^* \tilde {{\mathbf S}}_2 {\mathbf V}_2\!-\!{\mathbf D}_1^{1/2}{\mathbf U}_1^* \tilde {{\mathbf S}}_1{\mathbf U}_2{\mathbf D}_2^{1/2})\non
&\quad~& {\mathbf Q}^{-{1/2}}(\lambda_{p,j} {\mathbf I}\!-\!\tilde{{\mathbf F}})^{-1}{\mathbf Q}^{-{1/2}}(\lambda_{p,j} {\mathbf V}_2^*\tilde {{\mathbf S}}_2 {\mathbf V}_1\!-\!{\mathbf D}_2^{1/ 2}{\mathbf D}_2^* \tilde {{\mathbf S}}_1 {\mathbf U}_1 {\mathbf D}_1^{1/2})\non
\end{eqnarray}
with ${\mathbf Q}={\mathbf V}_2^* \tilde {{\mathbf S}}_2 {\mathbf V}_2$ and $ \tilde{{\mathbf F}}= {n_1}^{-1}{\mathbf Q}^{-{1/2}}{\mathbf D}_2^{1/2}{\mathbf U}_2^* {\mathbf X}{\mathbf X}^*{\mathbf U}_2{\mathbf D}_2^{1/2}{\mathbf Q}^{-{1/2}}$.
\begin{lem}\label{lem1}
Assume that ${\mathbf K}(\lambda_{p,j})$ is defined in (\ref{Kn}). Then,
\begin{align}
 {\mathbf K}(\lambda_{p,j})-\psi_k\um(\psi_k) {\mathbf D}_1-c_2\psi_k^2m(\psi_k){\mathbf I}_M-\psi_k {\mathbf I}_M \overset{a.s.}{\rightarrow} {\mathbf 0}_{M\times M}
\end{align}
where $\psi_k=:\psi(\alpha_k)$ is the limit of $\lambda_{p,j}$, $m(\cdot)$ is the Stieljtes transform of $ \tilde {\mathbf F}$ and $\um(\lambda)=-(1- c_1)/{\lambda}+ c_1 m(\lambda)$.
\end{lem}


According to Lemma \ref{lem1}, we obtain that $\psi_k$ satisfies the following equation:
\begin{equation}
\psi_{k}+c_{2}\psi^2_{k}m(\psi_{k})+\psi_k\um(\psi_{k}) \alpha_k=0.\label{0eqa}
\end{equation}
Therefore, the estimator of the population spiked eigenvalue, $\alpha_k$, is obtained as below:
\be
\hat \alpha_k= -\displaystyle\frac{1+ c_2 \lambda_{p,j} m(\lambda_{p,j})}{\um(\lambda_{p,j})},\label{ha1}
\ee
where $j \in \mathcal{J}_k, k=1,\cdots, K$ and  $m(\cdot)$   is approximately the same as
the Stieltjes transform of the LSD of  the Fisher matrix ${\mathbf F}$ if  the number of its spikes is fixed.

Next, the estimates of  $m(\lambda_{p,j})$ and  $\um(\lambda_{p,j})$ in (\ref{ha1}) are also provided.
We adopt an approach similar to that in \cite{JiangBai2019} to estimate $m(\lambda_{p,j})$.
Define $r_{ij}=|\lambda_{p,i}-\lambda_{p,j}|/|\lambda_{p,j}|$ and the set $\mathcal {J}_0=\{i\in (1,\cdots, p): r_{ij}\leq 0.2 \}$ and    $\tilde c_\ell=(p-|\mathcal {J}_0|)/n_\ell, \ell=1,2.$; then,
\be\hat m(\lambda_{p,j})=\frac{1}{p-|\mathcal {J}_0|}\sum\limits_{i\notin \mathcal J_0}(\lambda_{p,i}-\lambda_{p,j})^{-1}\label{procedure}
\ee
is a good estimator of $m(\lambda_{p,j})$,
where the set $\mathcal {J}_0$ is  selected to avoid  the effect of multiple roots and to make the estimator more accurate.
%
Furthermore,  the estimator of $\um(\lambda_{p,j})$  is obtained by the equation
\[\hat \um (\lambda_{p,j})=-\displaystyle\frac{1-\tilde c_1}{\lambda_{p,1}}+\tilde c_1 \hat m(\lambda_{p,j})\]
The estimator  $\hat \alpha_k$ is calculable in practice and is expressed as
\be
\hat \alpha_k= -\displaystyle\frac{1+\tilde c_2 \lambda_{p,j} \hat m(\lambda_{p,j})}{\hat \um(\lambda_{p,j})}.\label{hak}
\ee

\section{Simulation Study}\label{Sim}

We conduct simulations that support the theoretical results and illustrate the accuracy of the estimators of the population distant spiked eigenvalues. Assume $p=100,200,400$, $n_1=2p, n_2=4p$ and
the matrix ${\mathbf T}_p{\mathbf T}_p^*$ is a general positive definite matrix  satisfying
${\boldsymbol\Sigma}_2={\mathbf I}_p$ and ${\boldsymbol\Sigma}_1={\mathbf U}_0 {\boldsymbol\Lambda} {\mathbf U}_0^*$, where  $\boldsymbol\Lambda$ is a diagonal matrix with the form
\begin{align*}
10, 7.5,7.5, \underbrace{2,\cdots, 2}_{(p-6)/2}, \underbrace{1,\cdots, 1}_{(p-6)/2}, 0.2,0.2,0.1.
\end{align*}

Here, $\alpha_1=10$, $\alpha_2=7.5$, $\alpha_3=0.2$ and $\alpha_4=0.1$. Let ${\mathbf U}_0$  be equal to the matrix composed of eigenvectors
of the following matrix
\begin{align}
\left(
\begin{array}{ccccc}
1 & \rho & \rho^2 & \cdots & \rho^{p-1} \\
\rho & 1 & \rho & \cdots & \rho^{p-2} \\
\ldots & \ldots & \ldots & \ldots &\ldots  \\
\rho^{p-1} & \rho^{p-2} & \rho^{p-3} & \cdots & 1 \\
\end{array}
\right),
\end{align}
where $\rho=0.5$.
We propose that the samples are from three kinds of populations. In detail, $x_{ij}$ and $y_{ij}$  are the
$i.i.d.$ samples from the Gaussian distribution, the chi-square distribution and the uniform distribution with mean $0$ and variance $1$.
Then, the frequency histograms of the estimators ${\hat{\alpha}}_i $, $i\in{\left\lbrace 1,2,3,4\right\rbrace }$  are depicted in the following figures using 5000 repetitions.
\begin{figure}[htbp]
	\centering
	\includegraphics[width=0.325\textwidth]{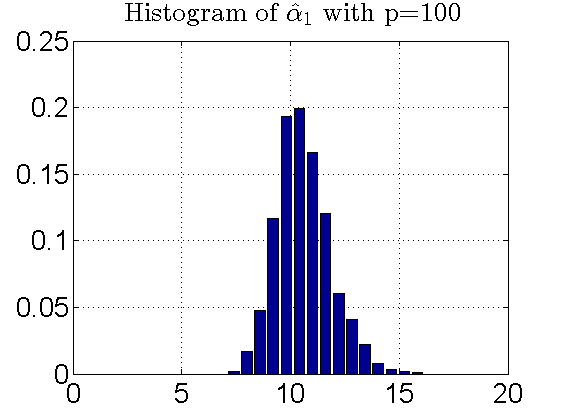}
	\includegraphics[width=0.325\textwidth]{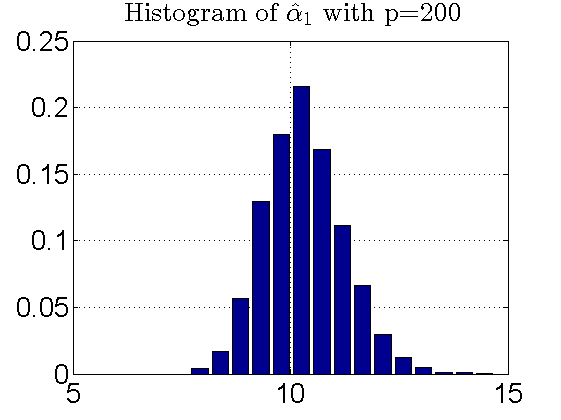}
	\includegraphics[width=0.325\textwidth]{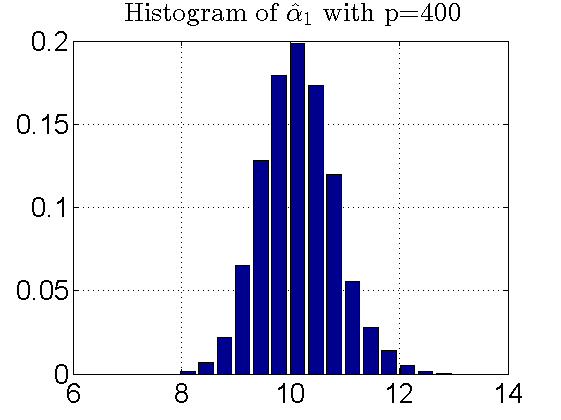}
	\caption{Estimating $\alpha_1$ $(\alpha_1=10)$ under the normal distribution assumption with $p=100, 200$ and $400$. }
\end{figure}\label{p1}

\begin{figure}[htbp]\label{p2}
	\centering
	\includegraphics[width=0.32\textwidth]{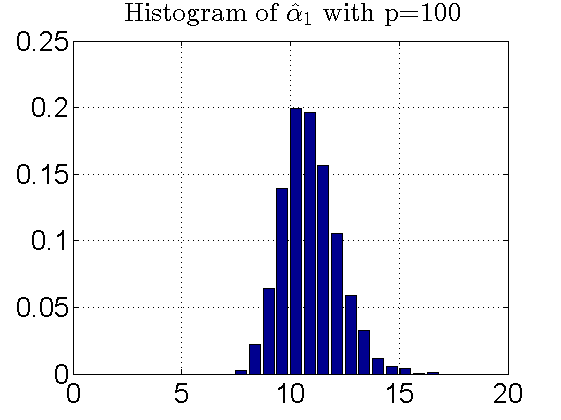}
	\includegraphics[width=0.32\textwidth]{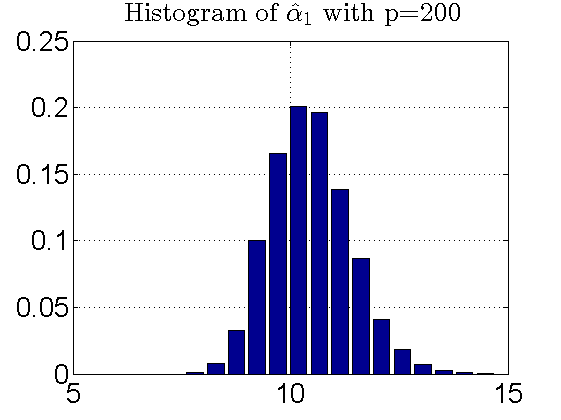}
	\includegraphics[width=0.32\textwidth]{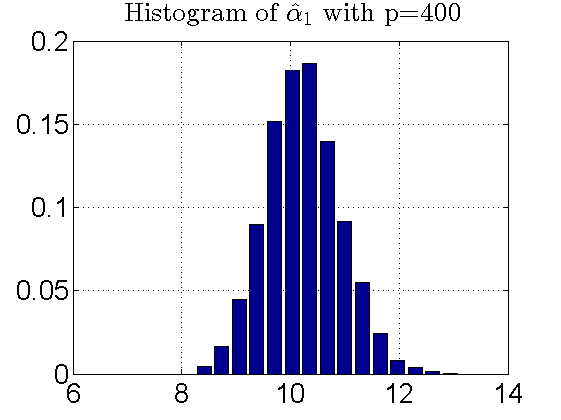}
	\caption{Estimating $\alpha_1$ $(\alpha_1=10)$ under the chi-square distribution assumption with $p=100, 200$ and $400$. }
\end{figure}

\begin{figure}[htbp]
	\centering
	\includegraphics[width=0.32\textwidth]{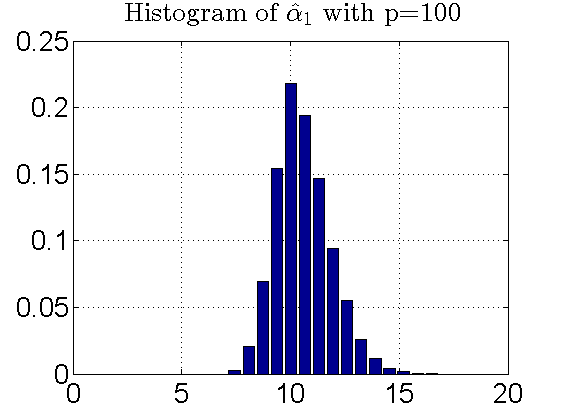}
	\includegraphics[width=0.32\textwidth]{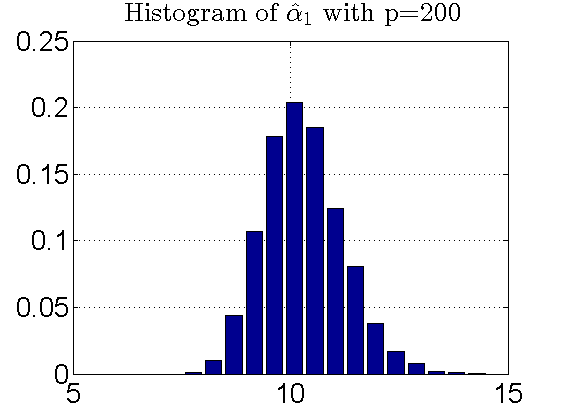}
	\includegraphics[width=0.32\textwidth]{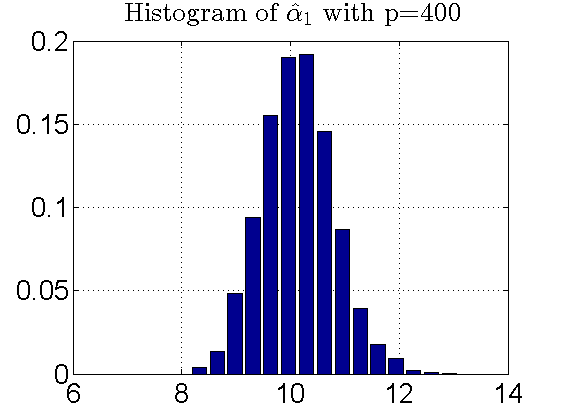}
	\caption{Estimating $\alpha_1$ $(\alpha_1=10)$ under the uniform distribution assumption with $p=100, 200$ and $400$. }
\end{figure}

\begin{figure}[htbp]
	\centering
	\includegraphics[width=0.32\textwidth]{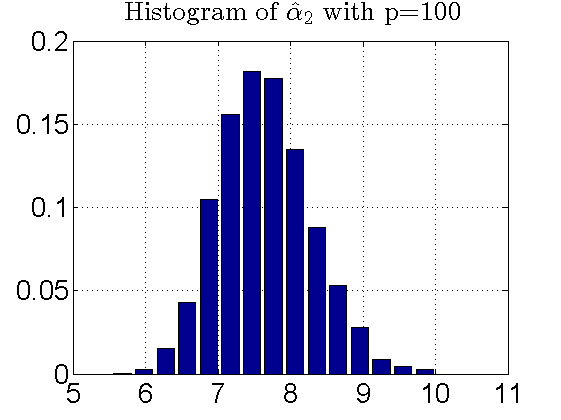}
	\includegraphics[width=0.32\textwidth]{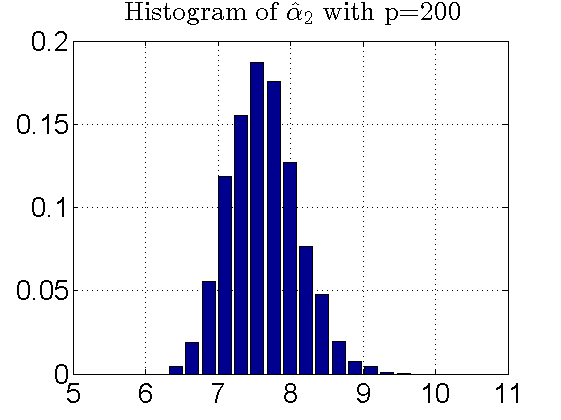}
	\includegraphics[width=0.32\textwidth]{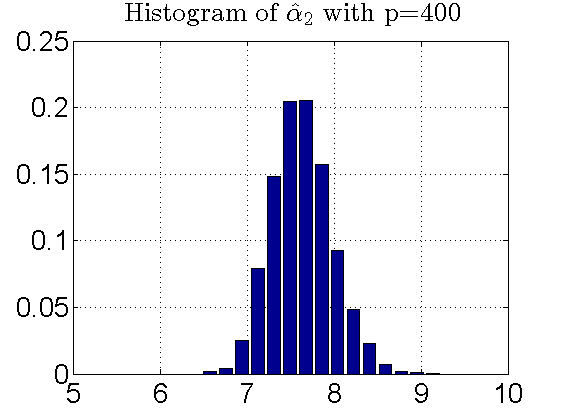}
	\caption{Estimating $\alpha_2$ $(\alpha_2=7.5)$ under the normal distribution assumption with $p=100$, $200$ and $400$. }
\end{figure}

\begin{figure}[htbp]
	\centering
	\includegraphics[width=0.32\textwidth]{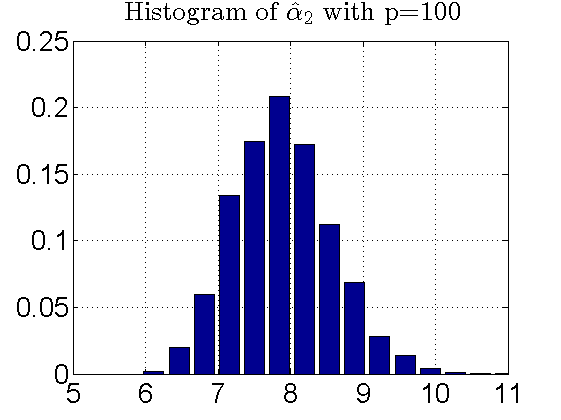}
	\includegraphics[width=0.32\textwidth]{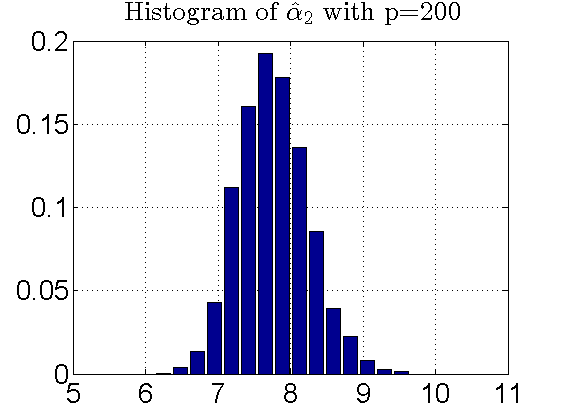}
	\includegraphics[width=0.32\textwidth]{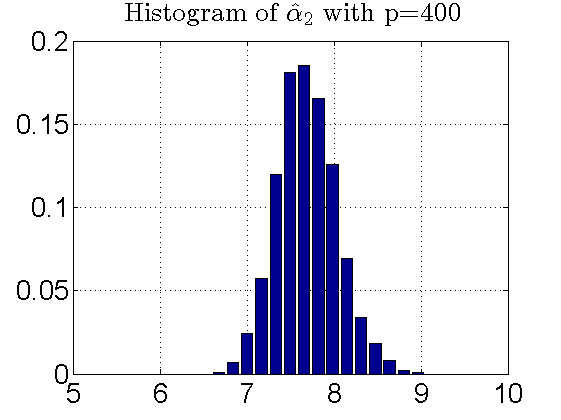}
	\caption{Estimating $\alpha_2$ $(\alpha_2=7.5)$ under the chi-square distribution assumption with $p=100, 200$ and $400$. }
\end{figure}

\begin{figure}[htbp]
	\centering
	\includegraphics[width=0.32\textwidth]{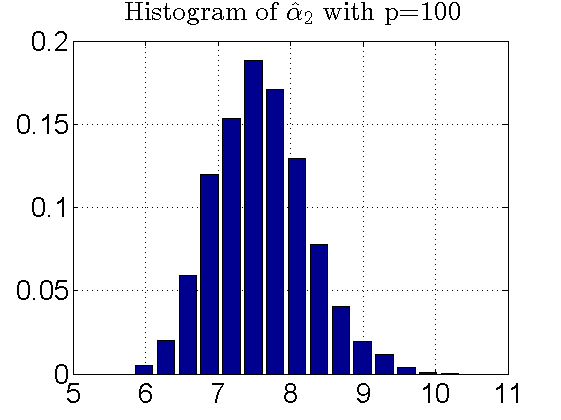}
	\includegraphics[width=0.32\textwidth]{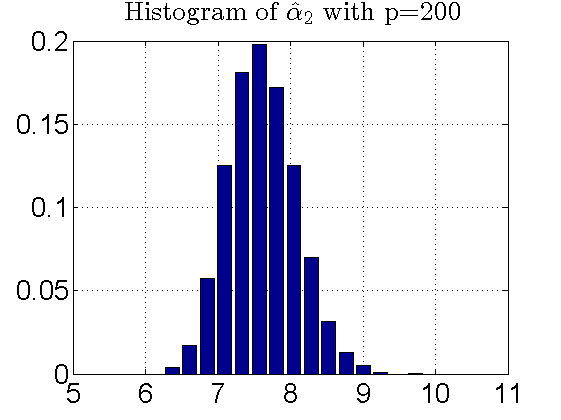}
	\includegraphics[width=0.32\textwidth]{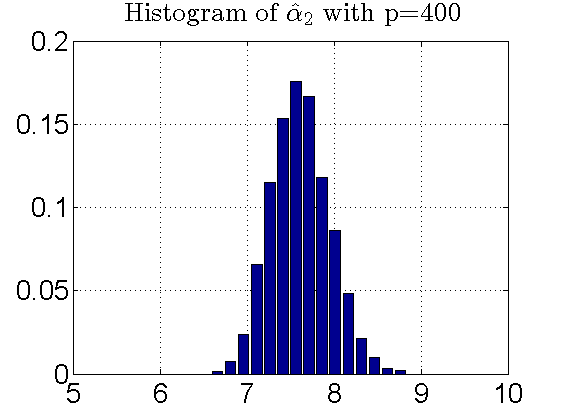}
	\caption{Estimating $\alpha_2$ $(\alpha_2=7.5)$ under the uniform distribution assumption with $p=100, 200$ and $400$. }
\end{figure}

\begin{figure}[htbp]
	\centering
	\includegraphics[width=0.32\textwidth]{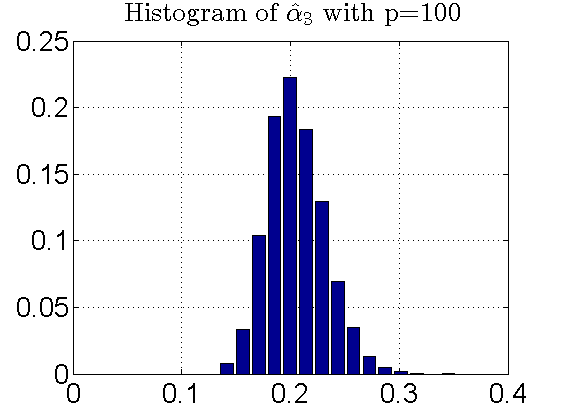}
	\includegraphics[width=0.32\textwidth]{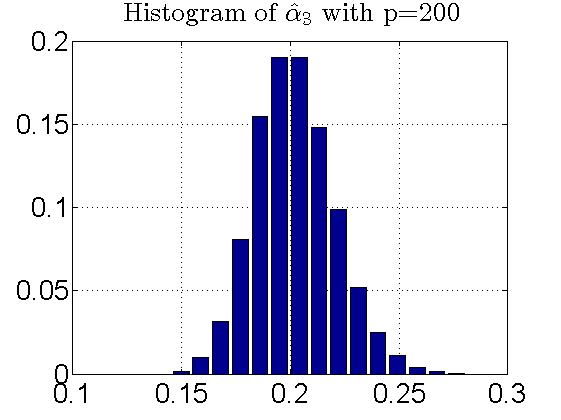}
	\includegraphics[width=0.32\textwidth]{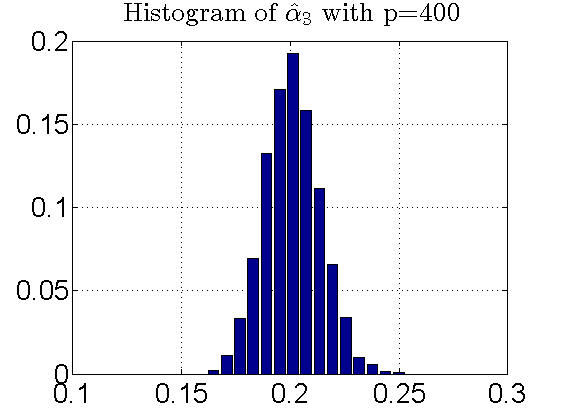}
	\caption{Estimating $\alpha_3$ $(\alpha_3=0.2)$ under the normal distribution assumption with $p=100, 200$ and $400$. }
\end{figure}

\begin{figure}[htbp]
	\centering
	\includegraphics[width=0.32\textwidth]{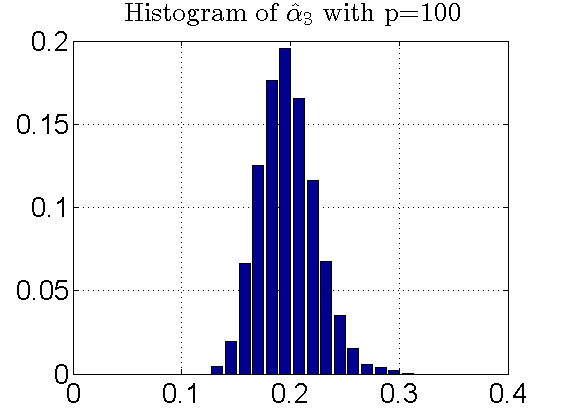}
	\includegraphics[width=0.32\textwidth]{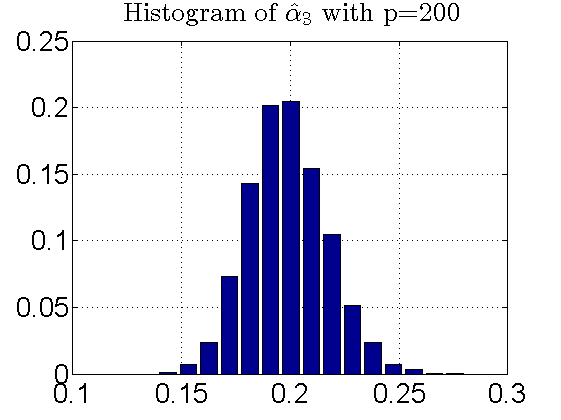}
	\includegraphics[width=0.32\textwidth]{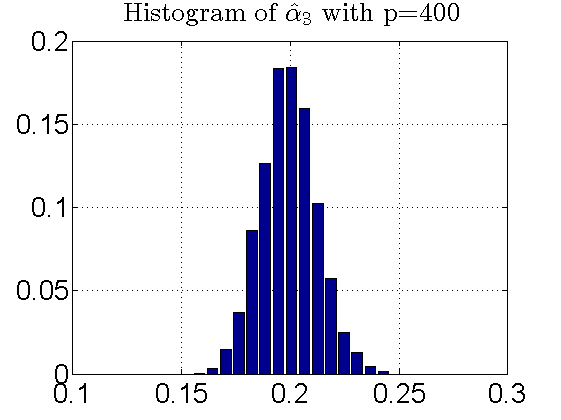}
	\caption{Estimating $\alpha_3$ $(\alpha_3=0.2)$ under the chi-square distribution assumption with $p=100, 200$ and $400$. }
\end{figure}

\begin{figure}[htbp]
	\centering
	\includegraphics[width=0.32\textwidth]{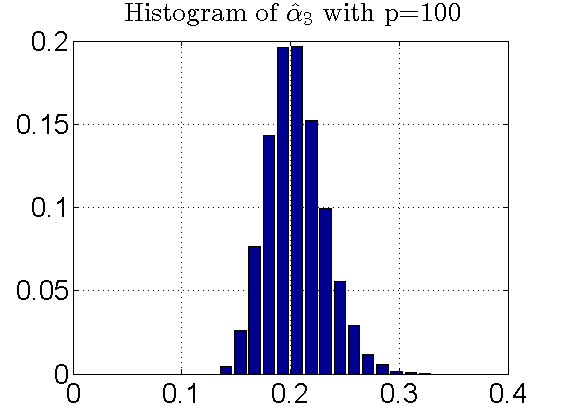}
	\includegraphics[width=0.32\textwidth]{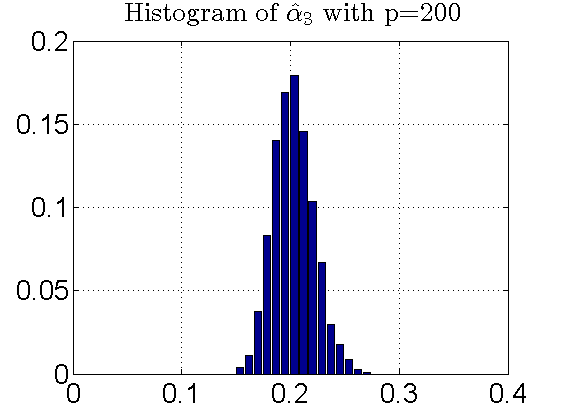}
	\includegraphics[width=0.32\textwidth]{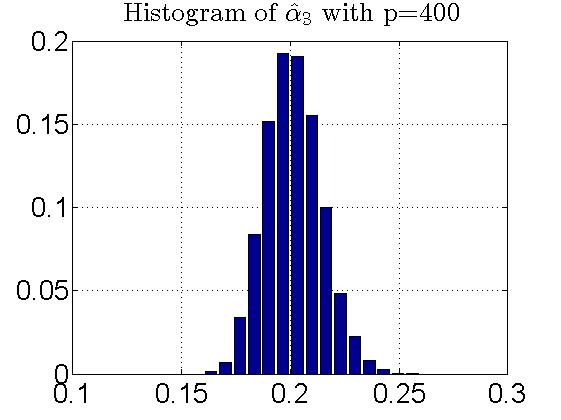}
	\caption{Estimating $\alpha_3$ $(\alpha_3=0.2)$ under the uniform distribution assumption with $p=100, 200$ and $400$. }
\end{figure}

\begin{figure}[htbp]
	\centering
	\includegraphics[width=0.3\textwidth]{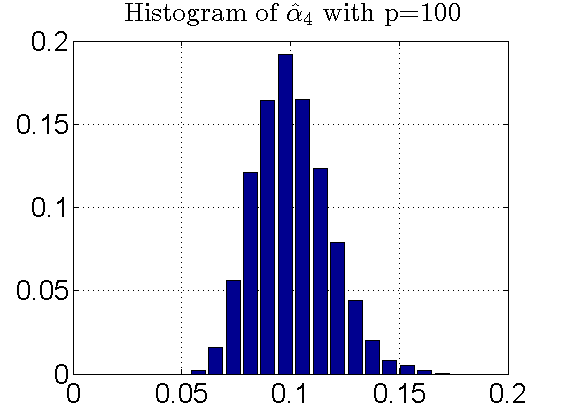}
	\includegraphics[width=0.3\textwidth]{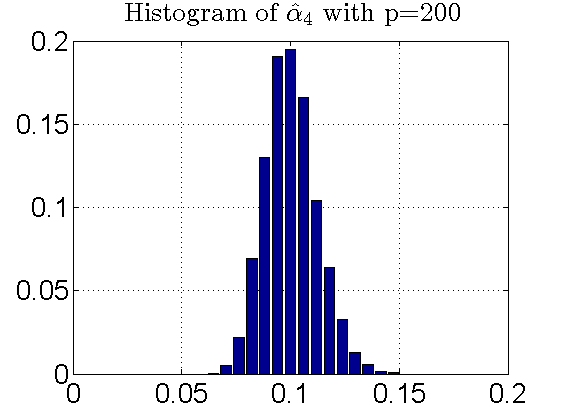}
	\includegraphics[width=0.3\textwidth]{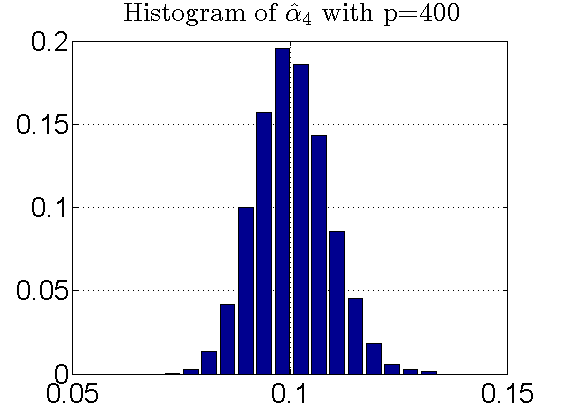}
	\caption{Estimating $\alpha_4$ $(\alpha_4=0.1)$ under the normal distribution assumption with $p=100, 200$ and $400$. }
\end{figure}

\begin{figure}[htbp]
	\centering
	\includegraphics[width=0.3\textwidth]{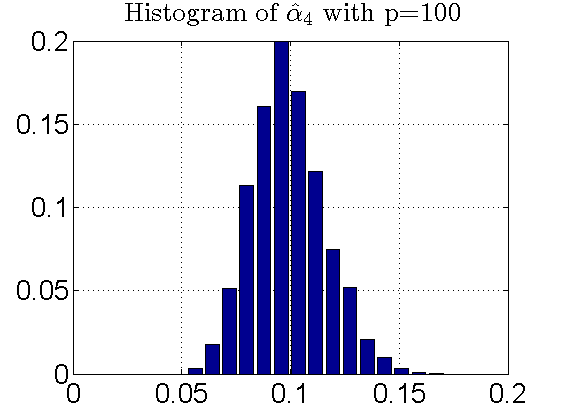}
	\includegraphics[width=0.3\textwidth]{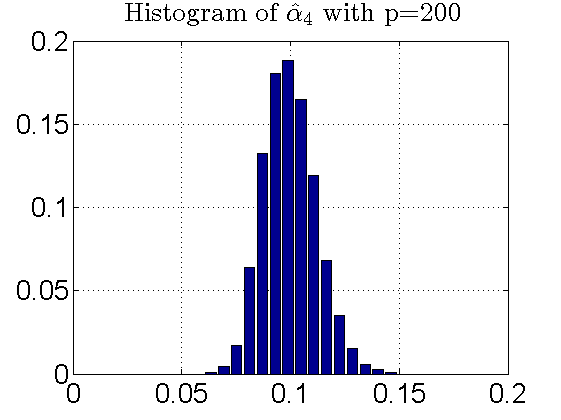}
	\includegraphics[width=0.3\textwidth]{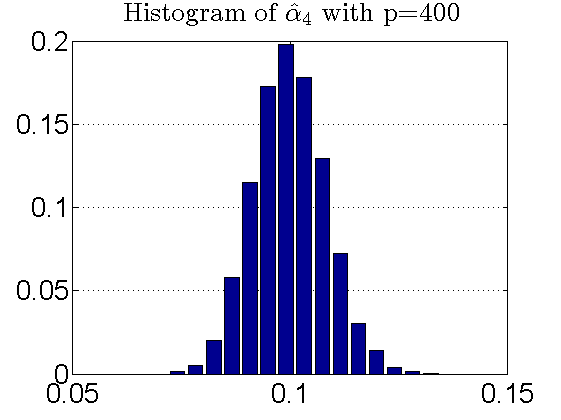}
	\caption{Estimating $\alpha_4$ $(\alpha_4=0.1)$ under the chi-square distribution assumption with $p=100, 200$ and $400$. }
\end{figure}

\begin{figure}[htbp]
	\centering
	\includegraphics[width=0.3\textwidth]{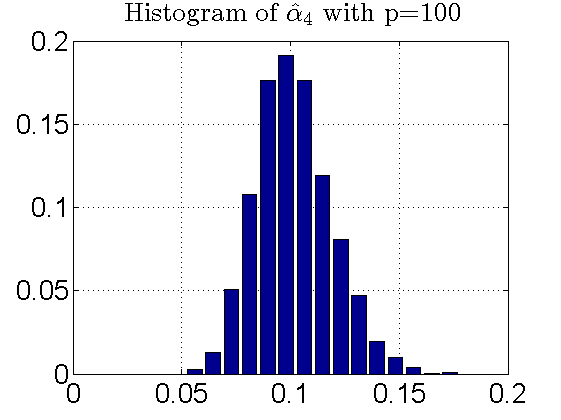}
	\includegraphics[width=0.3\textwidth]{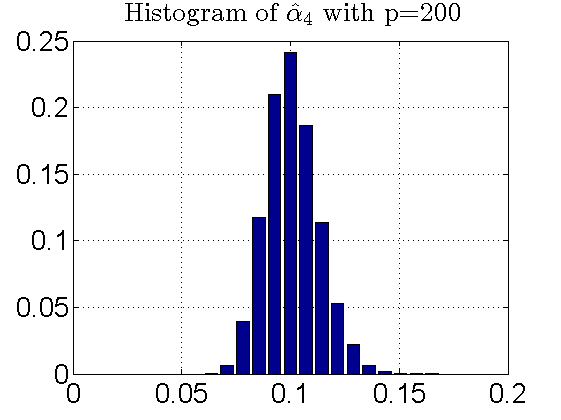}
	\includegraphics[width=0.3\textwidth]{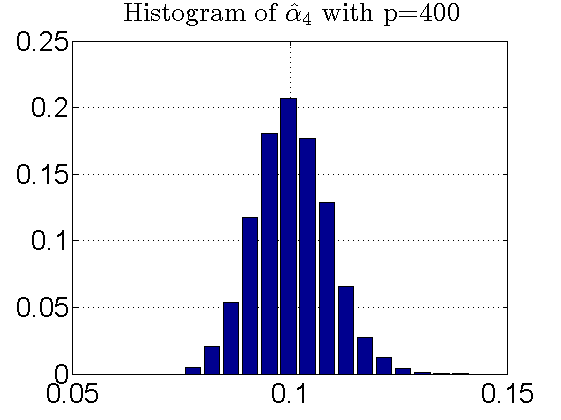}
	\caption{Estimating $\alpha_4$ $(\alpha_4=0.1)$ under the uniform distribution assumption with $p=100, 200$ and $400$.}
\end{figure}
Figures~1, 4, 7, and 10 show the accuracy of estimating $\alpha_i$, $i\in{\left\lbrace 1,2,3,4\right\rbrace }$ with $x_{ij}$ and $y_{ij}$ being drawn independently from $\mathcal{N}(0,1)$;
Figures~2, 5, 8, and 11 show the accuracy of estimating $\alpha_i$, $i\in{\left\lbrace 1,2,3,4\right\rbrace }$ with $x_{ij}$ and $y_{ij}$ being drawn independently  from $\chi^2(2)/2-1$; and
Figures~3, 6, 9, and 12 show the accuracy of estimating $\alpha_i$, $i\in{\left\lbrace 1,2,3,4\right\rbrace }$ with $x_{ij}$ and $y_{ij}$ being drawn independently  from  $U(-\sqrt{3},\sqrt{3})$.
For the single roots $\alpha_1$ and $\alpha_4$, the (\ref{hak}) are applied to the largest and the least sample eigenvalues, respectively.
For the multiple roots $\alpha_2$ and $\alpha_3$,
we first estimate the spike $\alpha_2$ with the second and third largest sample eigenvalues, respectively, and then take their average to obtain the final estimate of the corresponding spike.
The estimator of the spike $\alpha_3$ can be obtained  by the sample eigenvalues $\lambda_{p,p-2}$  and $\lambda_{p,p-1}$  in a similar way.
As seen from the figures, we find that the accuracy of estimates of the spikes  improves more and that the range of each  estimator decreases
as  the dimensionality $p$ increases  under all three distribution assumptions. In other words,  the estimates are
more focused and accurate when the  dimensionality $p$ continues to increase.

\section{Conclusion}\label{Con}
In this paper,  the phase transition of the spikes for a generalized Fisher matrix is proposed.
We extend the result in \cite{WangYao2017} to a general case to better match actual cases.
More importantly, the estimates of the population spiked eigenvalues are also provided, and thus, our results are calculable
and feasible in practice. As is known, the phase transition is the basis for the study of the asymptotic distribution for the  sample spiked eigenvalues.
In future work, we will investigate the CLT  in a high-dimensional Fisher matrix.

%
%
%

\section*{References}

\bibliography{mybibfile}

\newpage


\begin{center}
	{\bf Supplement to ``The limits of the distant  sample spikes for a high-dimensional generalized Fisher matrix and its applications".}
\end{center}
\renewcommand\thesection{\Alph{section}}
\setcounter{section}{0}


\section{Proof of Lemma \ref{lem1}}
Based on the expression of ${\mathbf K}(\lambda_{p,j})$ defined in (\ref{Kn}), we have
 \begin{eqnarray*}
{\mathbf K}(\lambda_{p,j})&=&\lambda_{p,j} {\mathbf V}_1^*\tilde {{\mathbf S}}_2 {\mathbf V}_1\!-\!{\mathbf D}_1^{1/2}{\mathbf U}_1^* \tilde {{\mathbf S}}_1{\mathbf U}_1{\mathbf D}_1^{1/2}\!-\!(\lambda_{p,j} {\mathbf V}_1^* \tilde {{\mathbf S}}_2 {\mathbf V}_2\!-\!{\mathbf D}_1^{1/2}{\mathbf U}_1^* \tilde {{\mathbf S}}_1{\mathbf U}_2{\mathbf D}_2^{1/2})\non
&\quad~& {\mathbf Q}^{-{1/2}}(\lambda_{p,j} {\mathbf I}\!-\!\tilde{{\mathbf F}})^{-1}{\mathbf Q}^{-{1/2}}(\lambda_{p,j} {\mathbf V}_2^*\tilde {{\mathbf S}}_2 {\mathbf V}_1\!-\!{\mathbf D}_2^{1/ 2}{\mathbf D}_2^* \tilde {{\mathbf S}}_1 {\mathbf U}_1 {\mathbf D}_1^{1/2})\non
&=&
\frac{\lambda_{p,j}}{n_2}{\mathbf V}_1^*{\mathbf Y}{\mathbf Y}^*{\mathbf V}_1\!-\!\frac{\lambda_{p,j}}{n_1} {\mathbf D}_1^{1 / 2}{\mathbf U}_1^*{\mathbf X}\big(\lambda_{p,j}{\mathbf I}_{n_1}\!-\!\underline{\tilde {\mathbf F}}\big)^{-1}{\mathbf X}^*{\mathbf U}_1{\mathbf D}_1^{1 / 2}\non
&&-\frac{\lambda_{p,j}^2}{n_2^2}{\mathbf V}_1^*{\mathbf Y}{\mathbf Y}^*{\mathbf V}_2{\mathbf Q}^{-{1/2}} (\lambda_{p,j} {\mathbf I}_{p-M}\!-\!\tilde {\mathbf F}\big)^{-1}{\mathbf Q}^{-{1/2}}{\mathbf V}_2^*{\mathbf Y}{\mathbf Y}^*{\mathbf V}_1\non
&&\!+\!\frac{\lambda_{p,j}}{n_2}{\mathbf V}_1^*{\mathbf Y}{\mathbf Y}^*{\mathbf V}_2{\mathbf Q}^{\!-\!{1/2}}\big(\lambda_{p,j} {\mathbf I}_{p-M}\!-\!\tilde {\mathbf F}\big)^{\!-\!1} {\mathbf Q}^{\!-\!{1/2}}
\frac{1}{n_1}{\mathbf D}_2^{1/2}{\mathbf U}_2^*{\mathbf X}{\mathbf X}^* {\mathbf U}_1{\mathbf D}_1^{1/2}\non
&&\!+\!\frac{\lambda_{p,j}}{n_1}{\mathbf D}_1^{1 / 2}{\mathbf U}_1^* {\mathbf X}{\mathbf X}^*{\mathbf U}_2{\mathbf U}_2^{1/2}{\mathbf Q}^{-{1\over2}} (\lambda_{p,j} {\mathbf I}_{p-M}\!-\!\tilde {\mathbf F})^{-1}{\mathbf Q}^{-{1\over 2}}\frac{1}{n_2} {\mathbf V}_2^*{\mathbf Y}{\mathbf Y}^* {\mathbf V}_1
\end{eqnarray*}
with ${\mathbf Q}={\mathbf V}_2^* \tilde {{\mathbf S}}_2 {\mathbf V}_2$ and $ \tilde{{\mathbf F}}= {n_1}^{-1}{\mathbf Q}^{-{1/2}}{\mathbf D}_2^{1/2}{\mathbf U}_2^* {\mathbf X}{\mathbf X}^*{\mathbf U}_2{\mathbf D}_2^{1/2}{\mathbf Q}^{-{1/2}}$.

According to the Fourth Moment Theorem in \cite{TaoVu2011}, the lemma 9.1 in \cite{BS2009} and the Borel-Cantelli lemma,  we can prove that the following convergence of matrices formula almost sure convergence.
The proof is mechanical and  tedious, and therefore, it is omitted here.
\begin{align}
&\frac{\lambda_{p,j}}{n_1 n_2}{\mathbf V}_1^*{\mathbf Y}{\mathbf Y}^*{\mathbf V}_2 {\mathbf Q}^{-{1 \over 2}}\big(\lambda_{p,j} {\mathbf I}_{p-M}- \tilde {\mathbf F}\big)^{-1} {\mathbf Q}^{-{1 \over 2}}
{\mathbf D}_2^{1 \over 2}{\mathbf U}_2^*{\mathbf X}{\mathbf X}^* {\mathbf U}_1 {\mathbf D}_1^{1 \over 2} \overset{a.s.}{\rightarrow} {\mathbf 0}_{M\times M}\label{eq01}\\
&\frac{\lambda_{p,j}}{n_1n_2}{\mathbf D}_1^{1 \over 2}{\mathbf U}_1^* {\mathbf X}{\mathbf X}^*{\mathbf U}_2{\mathbf D}_2^{1 \over 2}{\mathbf Q}^{-{1 \over 2}} (\lambda_{p,j} {\mathbf I}_{p-M}- \tilde {\mathbf F})^{-1}{\mathbf Q}^{-{1 \over 2}}{\mathbf V}_2^*{\mathbf Y}{\mathbf Y}^* {\mathbf V}_1
\overset{a.s.}{\rightarrow} {\mathbf 0}_{M\times M}\\ \label{eq02}
&\frac{\lambda_{p,j}}{n_2}{\mathbf V}_1^*{\mathbf Y}{\mathbf Y}^*{\mathbf V}_1\!-\!\psi_k {\mathbf I}_M \overset{a.s.}{\rightarrow} {\mathbf 0}_{M\times M}
\\
\!-\!&\frac{\lambda_{p,j}^2}{n_2^2}{\mathbf V}_1^*{\mathbf Y}{\mathbf Y}^*{\mathbf V}_2{\mathbf Q}^{\!-\!{1 \over 2}} (\lambda_{p,j} {\mathbf I}_{p-M}\!-\!\tilde {\mathbf F}\big)^{\!-\!1} {\mathbf Q}^{\!-\!{1 \over 2}}{\mathbf V}_2^*{\mathbf Y}{\mathbf Y}^*{\mathbf V}_1\!-\!c_2\psi_k^2m(\psi_k){\mathbf I}_M  \overset{a.s.}{\rightarrow} {\mathbf 0}_{M\times M}  \label{eq04}\\
-&\frac{\lambda_{p,j}}{n_1} {\mathbf D}_1^{1 \over 2}{\mathbf U}_1^*{\mathbf X}\big(\lambda_{p,j}{\mathbf I}_{n_1}\!-\!\underline{\tilde {\mathbf F}}\big)^{\!-\!1}{\mathbf X}^*{\mathbf U}_1{\mathbf D}_1^{1 \over 2}\!-\!\psi_k\um(\psi_k) {\mathbf D}_1 \overset{a.s.}{\rightarrow} {\mathbf 0}_{M\times M}
\label{eq05}
\end{align}

\section{Proof of Theorem \ref{F-limits}}
	
	For the generalized Fisher matrix ${\mathbf F}$  formulated in (\ref{F}),
	where $\tilde {{\mathbf S}}_1=n_1^{-1}{\mathbf X}{\mathbf X}^*$ and $\tilde {{\mathbf S}}_2={n_2}^{-1}{\mathbf Y}{\mathbf Y}^*$
	are the standardized sample covariance matrices. 
	Denote the Stieltjes transform of the LSD of the matrix ${\mathbf F}$ as $m_{{\mathbf c}, H}(\lambda)$ and that of matrix $\underline{{\mathbf F}}=n_1^{-1}{\mathbf X}^*{\mathbf T}_{p}^*\tilde{{\mathbf S}}_2^{-1}{{\mathbf T}}_{p}{\mathbf X}$
	as ${\underline m}_{{\mathbf c}, H}(\lambda)$. The LSD of $\tilde{{\mathbf S}}_2$ is presented as $F_{c_2}$  and its
	Stieltjes transform is $m_{c_2}(\lambda)$. Similarly, the LSD of $\tilde{{\mathbf S}}_2={n_2}^{-1}{\mathbf Y}^*{\mathbf Y}$
	is ${\underline F}_{c_2}$, which has the  Stieltjes transform denoted as $\um_{c_2}(\lambda)$.
	
	
	Furthermore,    for the nonzero spiked eigenvalues $\lambda_j \rightarrow \psi_k, j \in \mathcal{J}_k$,
	it follows from  equation (\ref{0eqa}) that
	\begin{align}\label{eigeneq3}
		1+{c}_{2}\psi_k m(\psi_k)+\um(\psi_k) \alpha_k=0.
	\end{align}
	By the relationship
	$\um(\psi_k) =-(1-c_1)/{\psi_k}+c_1 m(\psi_k)$,
	we obtain  that
	\bqa
	\um (\psi_k)=-\frac{h^2}{c_1\alpha_k +c_2\psi_k}, \label{umak}
	\eqa
	where $h^2=c_1+c_2-c_1c_2$. Furthermore, by (9.14.7) in \cite{BS2009}, we know that the $\um(\psi_k)$ satisfies the following equation
	\bqa
	\psi_k&=&-\frac{1}{\um(\psi_k)}+c_1\int \frac{1}{t+\um(\psi_k)} \md F_{c_2}(t)\non
	&=&-\frac{1}{\um(\psi_k)}+c_1\Big[\frac{1}{c_2}\um_{c_2}\big\{-\um(\psi_k)\big\}-\frac{1-c_2}{c_2 \um(\psi_k)}\Big]\non
	&=&-\frac{h^2}{c_2\um(\psi_k)}+\frac{c_1}{c_2}m_0.
	\label{lumm0}
	\eqa
	where $m_0=\um_{c_2}(-\um(\psi_k))$, combine (\ref{umak}) and (\ref{lumm0}), it follows that
	\be
	m_0(\psi_k)=-\alpha_k.\label{m0ak}
	\ee
	For each of the sample eigenvalues $\lambda_j, j\in {\mathcal J}_k, k=1, \cdots, K$ of the generalized Fisher matrix ${\mathbf F}$, apply (\ref{m0ak}) to equation (2.9) in \cite{Zhengetal2017}; then, it is obtained that
	\bqa
	\psi_k&=&\displaystyle\frac{m_0\Big(h^2+c_1\big(c_2\displaystyle\int\displaystyle\frac{m_0}{t+m_0}\md H(t)-1\big)\Big)}{c_2\big(c_2\displaystyle\int\displaystyle\frac{m_0}{t+m_0}\md H(t)-1\big)}\non
	&=&\displaystyle\frac{1-c_1\displaystyle\int\displaystyle\frac{t}{t+m_0}\md H(t)}{c_2\displaystyle\int\displaystyle\frac{1}{t+m_0}\md H(t) -\displaystyle\frac{1}{m_0}}\non
	&=&\displaystyle\frac{\alpha_k\big(1-c_1\displaystyle\int\displaystyle\frac{t}{t-\alpha_k}\md H(t) \big)}{1+c_2\displaystyle\int\displaystyle\frac{\alpha_k}{t-\alpha_k}\md H(t) }\label{eqA5}
	\eqa
	Combined with Theorem \ref{thm2} and \ref{thm3}, we prove that the limit of the sample eigenvalues $\lambda_j, j\in {\mathcal J}_k$ associated with the distant spike $\alpha_k$ is $\psi_k$. The limit of the sample eigenvalues associated with the closed spike is the border of the support of the LSD of the Fisher matrix ${\tilde{\mathbf S}}_1{{\mathbf T}}_{p}^*{\tilde{\mathbf S}}_2^{-1}{{\mathbf T}}_{p}$. For the proof details of the limit of sample eigenvalues associated with the closed spike, we refer the readers to Theorem 4.2 in \cite{BaiYao2012}.  The limit of the sample  closed spiked  eigenvalues can be derived in parallel according to their method. Now, the proof of Theorem \ref{F-limits} is completed.

\section{Proof of Theorem \ref{thm2}}
Since $x_0\in  \mathcal{G}^c_{F^{\mathbf c, H}}$, there exists $\delta>0$ such that $(x_0-\delta,x_0+\delta)\subset \mathcal{G}^c_{F^{{\mathbf c}, H}}$.  Write $z=x+iv$ with $x\in(x_0-\delta,x_0+\delta)$ and $v>0$. Then, by (2.9) in \cite{Zhengetal2017}, the  following equation
holds:
\bqa
z&=&\frac{h^2m_0(z)}{c_2\left(-1+c_2\displaystyle\int\frac{m_0(z)dH(t)}{t+m_0(z)}\right)}+\frac{c_1}{c_2}m_0(z)\nonumber\\
&=&\frac{m_0(z)\left(1-c_1\displaystyle\int\frac{tdH(t)}{t+m_0(z)}\right)}{-c_2\displaystyle\int\frac{tdH(t)}{t+m_0(z)}-1+c_2}\non
&:=&\psi(-m_0(z)),\label{eq1}
\eqa
where $m_{0}(z)=\underline{m}_{c_2}(-\underline{m}_{{\mathbf c}, H}(z))$ and $\underline{m}_{c_2}(z)$ is the unique solution, with the same sign of the imaginary parts as that of $z$, to the following equation:
\bqa
z=-\frac1{\underline{m}_{c_2}(z)}+c_2\int\frac{1}{t+\underline{m}_{c_2}(z)}dH(t).\label{eq2}
\eqa
Additionally, by the definition of $m_0(z)$, we have
\bqa
\underline {m}_{{\mathbf c}, H}(z)=\frac1{m_0(z)}-c_2\int\frac{1}{t+m_0(z)}dH(t).\label{eq3}
\eqa
Write $m_0(z)=m_{01}+im_{02}$ and $\underline {m}_{{\mathbf c}, H}(z)=g_1(z)+ig_2(z)$; we have that $g_2(z)\to g_2(x_0)=0$ as $v\to 0$. By equation (\ref{eq3}), we have
\bqa
g_2(z)=-m_{02}(z)\Big(\frac1{|m_0(z)|^2}-c_2\int\frac{1}{|t+m_0(z)|^2}dH(t)\Big).\label{eq4}
\eqa
By the definition of $m_0(z)$,
$$
m_0(z)=\underline{m}_{c_2}(-\underline{m}_{{\mathbf c}, H}(z))=\int\frac{d{\underline F}_{c_2}(\lambda)}{\lambda+\underline{m}_{{\mathbf c}, H}(z)},
$$
where ${\underline F}_{c_2}$ is the limiting spectral distribution corresponding to the Stieltjes transform $\underline m_{c_2}$.

Thus,
$$
m_{02}(z)=-g_2(z)\int\frac{d\underline{F}_{c_2}(\lambda)}{|\lambda+\underline{m}_{{\mathbf c}, H}(z)|^2}.
$$
Therefore, as $v\to 0$, $m_{02}(z)\to 0$ and
$$
-m_{02}(z)/g_2(z)\to \int\frac{d\underline{F}_{c_2}(\lambda)}{(\lambda+\underline{m}_{{\mathbf c}, H}(x))^2}>0.
$$
On the other hand, by equation (\ref{eq4}), the same limit shows that
\bqa\label{eq5}
\Big(\frac1{m_0^2(x)}-c_2\int\frac{1}{(t+m_0(x))^2}dH(t)\Big)>0.
\eqa
This shows (i) and (ii)  with  $u_0=-m_0(x_0)$.

Taking the imaginary parts of both sides of the equation (\ref{eq1}), we have
\[
v\!=\!-m_{02}(z)\Big(\frac{1\!-\!c_2\!-\!(c_1(c_2\!-\!1)\!+\!c_2)\!\int\!\frac{t^2dH(t)}{|t\!+\!m_0(z)|^2}\!-\!c_1c_2\left|\!\int\!\frac{tdH(t)}{t\!+\!m_0(z)}\right|^2\!+\!2c_2m_{01}(z)\!\int\!\frac{tdH(t)}{|t\!+\! m_0(z)|^2}}{|-c_2\!\int\!\frac{tdH(t)}{t\!+\!m_0(z)}\!-\!1\!+\!c_2|^2}\Big)
\]

Dividing both sides of the above equation  by $-m_{02}$ and then letting $v\to0$, the right-hand side of the above tends to
\begin{align}\label{eq6}
 \Big(\frac{1\!-\!c_2\!-\!(c_1(c_2\!-\!1)\!+\!c_2)\int\frac{t^2dH(t)}{(t\!+\!m_0(x))^2}\!-\!c_1c_2\left(\int\frac{tdH(t)}{t\!+\!m_0(x)}\right)^2\!+\!2c_2m_{0}(x)\!\int\!\frac{tdH(t)}{(t\!+\! m_0(x))^2}}{(\!-\!c_2\!\int\!\frac{tdH(t)}{t\!+\!m_0(x)}\!-\!1\!+\!c_2)^2}\Big)>0.
\end{align}

Write $u_0=-m_0(x_0)$. Then,
\bqn
\psi'(u_0)&=&\Re(\psi'(u_0))\non
&=&\lim_{m_{02}(x_0+iv)\to0}\frac{\psi(-m_0(z))}{d(-m_0(z))}\Big|_{z=x_0+iv}\\
&=&\Re\Big(\lim_{m_{02}(z)\to 0}\frac{\psi(-m_{0}(z))-\psi(-m_{01}(z))}{-im_{02}(z)}\Big)\Big|_{z=x_0+iv}\\
&=&\lim_{m_{02}(z)\to 0}\Im\Big(\frac{\psi(-m_{0}(z))-\psi(-m_{01}(z))}{-m_{02}(z)}\Big)\Big|_{z=x_0+iv}\\
&=&\lim_{m_{02}(z)\to 0}\Big(\frac{\Im(\psi(-m_{0}(z))}{-m_{02}(z)}\Big)\Big|_{z=x_0+iv}
\eqn
Note that the right-hand side of the above equation is the same as that of (\ref{eq6}) and hence is positive.

\section{proof of Theorem \ref{thm3}}

By (i) - (iii), there exists a constant $\delta>0$ such that $(u_0-\delta,u_0+\delta)\subset \mathcal{G}^c_{H}$ and the conditions $(i)-(iii)$ hold for all $u\in(u_0-\delta,u_0+\delta)$.
For $u\in (u_0-\delta, u_0+\delta)$ and $w\in R$, set $-m=u+iw$; then, by condition $(i)$, $z=\psi(-m)$ is an analytic function in the space of $m$. Additionally, by $(iii)$, there is a unique inverse function $m=m(\psi)$ of $\psi(u)$ such that  $\psi=\psi(-m(\psi))$ for all $\psi \in(x_0-\eta,x_0+\eta)$ where $\eta>0$ is a constant. Since $\psi(-m)$ is analytic, its inverse function
is also analytic when $\psi\in(x_0-\eta,x_0+\eta)$; therefore, the inverse function can be extended to an open region $\mathcal{B}$ containing $(x_0-\eta,x_0+\eta)$ as a subset.

On the other hand, for all $z\in C^+$, by \cite{Zhengetal2017}, the Stieltjes transform $\underline m_{{\mathbf c}, H}(z)$ of the LSD of the Fisher matrix is uniquely determined by equations (\ref{eq1}) and (\ref{eq3}). Specifically, when $z\in \mathcal{B}$, $m_0(z)=m(z)$. When $v\to 0$, $\Im(m)\to 0$. Then, by (\ref{eq3}), $\Im(\underline m_{{\mathbf c}, H}(z))\to 0$, for all $\Re(z)\in(x_0-\eta, x_0+\eta)$. Therefore, $x_0\in \mathcal{G}^c_{F^{\mathbf c, H}}$.
\end{document}